\documentclass[11pt]{amsart}
\usepackage{amssymb, amscd, amsmath}

\newtheorem{thm}{Theorem}[section]
\theoremstyle{definition}
\newtheorem{defn}{Definition}[section]
\theoremstyle{remark}

\theoremstyle{plain}
\newtheorem{lem}[thm]{Lemma}
\newtheorem{cor}[thm]{Corollary}

\newtheorem{rem}[thm]{Remark}

\def\C{{\mathbb C}}                               
\def\Q{{\mathbb Q}}                               
\def\Z{{\mathbb Z}}                               
\def\O{{\mathcal O}}                              
\def\disc{\operatorname{disc}}                    
\def\det{\operatorname{det}}                      
\def\trace{\operatorname{trace}}                  

\def\t{{\times}}                                  
\def\s1{\sum_{i = 1}^{j_{\alpha}}({t_i}^3~-~ 3{t_i}{\Delta_i})}
\def\u2{\sum_{m = 1}^{l} {\tilde t_m}^3}
\def\v3{\sum_{j = 1}^{s} 3~{\tilde t'_j}}
\def\w4{\sum_{i=1}^{j_{\alpha}}({t_i}^4~-~4~{t_i}^2~{\Delta_i}~+~2~{\Delta_i^2})}
\def\x5{\sum_{m = 1}^{l} {\tilde t_m}^4}
\def\y6{\sum_{j = 1}^{s} 4~{\tilde {t'_j}}^2{\tilde {\Delta'_j}}}
\def\z7{\sum_{n = 1}^{r} 2~{\tilde {\Delta''_n}}^2}

\def\vpmod{\!\!\!\pmod}

\begin{document}
\title[Waring's problem for matrices]{Waring's problem for matrices over orders\\ in algebraic number fields}
\author[A. S. Gadre and S. A. Katre]{A. S. Gadre and S. A. Katre}
\address{Department of Mathematics, University of Pune, Ganeshkhind, Pune - 411 007, INDIA.}
\email[A. S. Gadre, S. A. Katre]
{\vskip1mm 
ganura@math.unipune.ernet.in, 
sakatre@math.unipune.ernet.in}
\date{\today}

\begin{abstract}
In this paper we give necessary and sufficient trace conditions for an $n \t n$ matrix over any commutative and
associative ring with unity to be a sum of $k$-th powers of matrices over that ring, where $n,k \geq 2$ are integers.
We prove a discriminant criterion for every $2 \t 2$ matrix over an order
$R$ in an algebraic number field to be a sum of cubes and fourth
powers of matrices over $R$. We also show that if $q$ is a prime
and $n \geq 2$, then every $n \times n$ matrix
over the ring $\O$ of integers of a {\it quadratic} number field 
is a sum of $q$-th powers (of matrices) over $\O$  
if and only if $q$ is coprime to the discriminant of the 
quadratic field.

\end{abstract}

\maketitle
\section{Notations.}
\begin{center}
 \begin{tabular}{rcl}
$R$ & = & A commutative and associative ring with unity (or \\
&&in particular an order in an algebraic number field).\\
$M_n(R)$ & = & The set of all $n \t n$ matrices over $R$.\\
$\Z$     & = & The set of integers.\\
$K$       & = & An algebraic number field.\\
$\O$      & = & The ring of integers of an algebraic number field $K$.\\
$\disc K$ & = & The discriminant of the algebraic number field $K$. \\
$\O_{m}$ & = & The additive subgroup of $\O$ generated by the $m$-th powers\\
                      & & of elements of $\O$, where $m$ is a positive integer.
 \end{tabular}
\end{center}

\section{Introduction.}
For us a ring will always mean a commutative and associative ring with
unity. In \cite{R} Richman proved the following theorem:

\begin{thm}\label{thm:Richman}
Let $R$ be a ring. If $n \geq k \geq 2$ are integers, then
the following statements are equivalent:
\begin{enumerate}
\item $M$ is a sum of $k$-th powers in $M_n(R)$;
\item $M$ is a sum of seven $k$-th powers in $M_n(R)$;
\item $M$ belongs to $M_n(R)$ and for every prime power $p^e$
dividing $k$, there are elements\\
$x_0=x_0(p), \dots, x_e = x_e(p)$ in $R$ such that
\begin{equation*}
{\trace}~ M = {x_0}^{p^e}+p{x_1}^{p^{e-1}}+ \cdots +p^{e}x_e.
\end{equation*}
\end{enumerate}
Moreover if $k = p$ is a prime, in condition $(2)$  
$``seven"$ can be replaced by $``five"$.  
Also condition $(3)$ simplifies to 
$``\trace M = {x_0}^{p}+ p~{x_1}"$, for some $x_0, x_1 \in R$.
\end{thm}

In \cite{N} M. Newman proved that every integral $2 \t 2$ matrix
is a sum of at most $3$ integral squares, and this
is the best possible. For every integer $n > 2$, 
he showed that every $M$ in $M_n(\Z)$ is a sum of squares of at most 
$7$ or $9$ integral
matrices, when $n$ is even and odd respectively.
Here he queried whether this problem could be generalised to the ring of 
integers of an algebraic number field. He also asked if the question could be generalised to the case of arbitrary powers.   

Motivated by these questions Richman in \cite{R} proved the following 
result: Let $n$ and $k$ be positive integers with $n,k \geq 2$. 
Then there exists a number $G_k$ (independent of $n$) 
such that every $n \t n$ matrix with integer entries is a sum of 
$G_k$ $k$-th powers in $M_n(\Z)$.

The Waring's problem for matrices over the ring $\O$ of integers 
of an algebraic number field $K$, was first treated by S. A. Katre and 
S. A. Khule in \cite{KK1}, where for $n \geq 2$,
they gave a discriminant criterion for $n \t n$ matrices over orders in algebraic number fields to be sums of squares. Later they also obtained (see \cite{KK2}) such a criterion for $k$-th powers $(n \geq k \geq 2)$. 

\begin{thm}\label{thm:ngek}(S. A. Katre-S. A. Khule, \cite{KK2}, Theorem $1$)
Let $R$ be an order in an algebraic number field $K$. Let 
$n \geq k \geq 2$. Then every $n \t n$ matrix over $R$ is a 
sum of (seven) $k$-th powers if and only if $(k, \disc R) =1$. 
\end{thm}

Note that the main results in \cite{KK1}, \cite{KK2} and \cite{R} 
make the assumption that $n \geq k$, that is the size of the matrix is greater than or equal to the power occurring in the decomposition.
This is because the proofs depend heavily on the construction 
of (an invertible) matrix $A$ in $M_n(R)$ whose characteristic polynomial is of the form ${(X-r)}^kX^{n-k}~-~1$, for some $r$ in $R$, with the
tacit assumption $n \geq k$. For $n < k$, this approach does not work. 
Thus arise the following natural questions:
\begin{enumerate}
\item If $n \geq 2$, (more specifically for $2 \leq n < k$),
what are the conditions under which an $n \t n$
matrix over a ring $R$ is a sum of $k$-th powers of matrices over $R$?  
\item If $2 \leq n < k$, what are the conditions under which 
every $n \t n$ matrix over an order $R$ in an algebraic number field
is a sum of $k$-th powers of matrices over $R$?
\end{enumerate}

\par In the light of the first question, in the next section,
for $n \geq 2$ we obtain trace conditions for a matrix $M$ in
$M_n(R)$ to be a sum of $k$-th powers in $M_n(R)$.
If a matrix $A \in M_n(R)$ is a sum of $k$-th powers of matrices
in $M_n(R)$ then $\trace A$ is a sum of traces of $k$-th
powers of matrices. In Theorem \ref{thm:Equiv} we show that
the converse of this also holds.
Note that Richman's conditions (see Theorem \ref{thm:Richman})
are simpler but are valid only for $n \geq k \geq 2$ whereas
Theorem \ref{thm:Equiv} proved here
is valid for $k \geq n \geq 2$ too.

In \cite{R}, Richman proved that for integers $n,k \geq 2$,
every $n \t n$ matrix over $\Z$ is a sum of $k$-th powers in $M_n{\Z}$.
Since $\Z$ is the ring of integers of the algebraic number field $\Q$,
one can ask if every $2 \t 2$ matrix over the ring $\O$ of integers of an algebraic number field $K$, can be written as a sum of $k$-th powers where $k \geq 3$. In the context of the ring of integers of an algebraic number field and orders therein we prove in sections
\S 5, \S 7, \S 8
the following results in connection with the second question:
\begin{enumerate}
\item Let $R$ be an order in an algebraic number field $K$. Then every $2 \t 2$
matrix over $R$ is a sum of cubes of matrices over $R$ if and only if
$(3, \disc R) = 1$.

\item Let $R$ be an order in an algebraic number field $K$. Then every $2 \t 2$ matrix over $R$ is a sum of fourth powers of matrices over $R$ 
if and only if $(2, \disc R) = 1$.

\item Let $n \geq 2$, $q \geq 3$ be a prime and $n < q$.
If $\O$ denotes the ring of integers of a {\it quadratic} number field 
then every $n \t n$ matrix over $\O$ 
is a sum of $q$-th powers of matrices over $\O$ 
if and only if $(q, \disc K) = 1$.
\end{enumerate}

\section{Trace conditions for matrices to be sums of $k$-th powers}

We recall a result below which highlights the significant role 
played by the trace of a matrix, while trying to decompose a given 
matrix as a sum of $k$-th powers of matrices.

\begin{thm}\label{thm:Trace}(D. R. Richman, \cite{R}, Proposition $5.5$)
Let $R$ be a ring. Let $n \geq 2$. Then, there exists a number 
$g_k$ (independent of $n$) such that every matrix in $M_n(R)$, 
whose trace lies in $k!R$ is a sum of $g_k$ $k$-th powers in $M_n(R)$.
\end{thm}

We derive some useful corollaries which will be used in proving the 
trace condition and other theorems in this paper.   

\begin{cor}\label{cor:T1}
Let $R$ be a ring.
Let $M$ belong to $M_n{R}$.
If $\trace~ M$ can be written as:
$\trace~ M =  \sum_{i = 1}^{l} \trace {M_i}^k + k!~r$, for
some $M_i \in M_n(R)$ and some $r \in R$, then 
$M$ is a sum of $k$-th powers of matrices in $M_n{R}$. 
\end{cor}

\begin{cor}\label{cor:T2}
Let $R$ be ring. Let $M$ belong to $M_n{R}$.
If every element of the quotient ring $R/{k!R}$ can be expressed 
as a sum of $k$-th powers of elements in $R/{k!R}$,
then every matrix $M$ in $M_n(R)$ is a sum of $k$-th powers of matrices in $M_n{R}$. 
\end{cor}
\begin{proof}
For $\alpha \in R$, ${\alpha}^k = \trace
((\rm{diag}(\alpha,0,\ldots,0))^{k})$.
\end{proof}

In a similar vein, we state a weaker version of the above result.

\begin{cor}\label{cor:T3}
Let $R$ be a ring and let $n,k$ be integers such that $n \geq 1, k \geq 2$. 
Let $M$ belong to $M_n{R}$.
If every element of $R$ can be written
as a sum of $k$-th powers of elements in $R$, 
then every matrix in $M_n(R)$ is a sum of $k$-th powers of matrices in $M_n{R}$. 
\end{cor}

The arithmetic question that arises in this context is to determine 
the condition when every element of an arbitrary ring $R$ 
can be written as a sum of traces of $k$-th powers of matrices in $M_n(R)$.
We investigate more about this in the lemma below.
 
\begin{lem}\label{lem:RGroup}
Let $R$ be a ring and $n, k \geq 2$ be integers.
Then, the set $$T = \bigg \{ \alpha \in R~|~ \alpha \equiv \sum_{j = 1}^{l} {\rm trace} (M_j)^k \vpmod{k!R},~{\rm for~ some~} M_1, \ldots, M_l ~{\rm in}~ M_n(R) \bigg \}$$
is an abelian group under addition.
\end{lem}

\begin{proof}
If $\alpha \in T$, then $-\alpha \equiv (k! -1)\alpha~\vpmod {k!R}$,
so $-\alpha$ belongs to $T$.  
\end{proof}

\begin{thm}\label{thm:Equiv}
Let $R$ be a ring and $n, k \geq 2$ be integers. 
Let $M$ belong to $M_n(R)$. Then the following are equivalent:
\begin{enumerate}
\item $M$ is a sum of $k$-th powers of matrices in $M_n(R)$.
\item $\trace~ M$ is a sum of traces of $k$-th powers 
of matrices in $M_n(R)$.
\item $\trace~ M$ is in the subgroup of $R$ generated by the traces of 
$k$-th powers of matrices in $M_n(R)$.
\item $\trace~ M$ is in the subgroup of $R$ generated by the
traces of $k$-th powers of matrices in $M_n(R)$ $\vpmod{k!R}$.
\item $\trace~ M$ is a sum of traces of $k$-th powers 
of matrices in $M_n(R)$ $\vpmod{k!R}$.
\end{enumerate}
\end{thm}

\begin{proof}
The implications $(1) \Longrightarrow (2) \Longrightarrow (3) \Longrightarrow (4)$ are clear.

$(5) \Longrightarrow (1)$. 
Follows from Corollary \ref{cor:T1}.
\end{proof}

\begin{cor}\label{cor:Fequiv}
Let $R$ be a ring and $n,k \geq 2$ be integers.
Then the following are equivalent:
\begin{enumerate}
\item Every matrix $M$ in $M_n(R)$ is a sum of $k$-th powers of matrices 
in $M_n(R)$.
\item Every element of $R$ is a sum of traces of $k$-th powers 
of matrices in $M_n(R)$.
\item $R$ is generated as a group by the
traces of $k$-th powers of matrices in $M_n(R)$.
\item $R$ is generated as a group by the
traces of $k$-th powers of matrices in $M_n(R)$ $\vpmod {k!R}$.
\item Every element of $R$ is a sum of traces of $k$-th powers 
of matrices in $M_n(R)$ $\vpmod {k!R}$.
\end{enumerate}
\end{cor}

\begin{proof}
The proof follows from Theorem \ref{thm:Equiv}, as given an  
$\alpha \in R$, $\alpha = \trace M$, for $M = 
\rm{diag}(\alpha, 0, \ldots, 0)$.
\end{proof}

\begin{lem}\label{lem:Ringqth}
Let $n,q \geq 2$ be integers and let $q$ be a prime.
If $R$ is a ring and $M$ belongs to $M_n(R)$, then 
$\trace (M^q) \equiv (\trace M)^q \vpmod{qR}$.
Also if every $n \t n$ matrix over $R$ is a sum of $q$-th
powers of matrices over $R$, then every element
of $R$ is a $q$-th power $\vpmod{qR}$.
\end{lem}
\begin{proof}
The proof follows from the part $(1) \Longrightarrow (3)$
in the proof of Proposition $2.5$, \cite{R} and 
that for $\alpha \in R$, $\alpha = \trace(\rm{diag}(\alpha, 0, \ldots, 0))$. 
Note that the condition $``n \geq p"$ in the assumption of 
Proposition $2.5$ of \cite{R} is not used in the proof 
$(1) \Longrightarrow (3)$. 
\end{proof}

\section{Discriminant conditions for orders}
We recall in this section a few results from \cite{KK2}. 
First we state the ``Discriminant criterion'' which is a 
crucial input for all the theorems proved in the later 
sections of the paper.  

\begin{lem}\label{lem:Disc}(\cite{KK2}, Lemma $3$)
Let $p$ in $\Z$ be a prime number. 
Let $\O$ be the ring of integers of a number field $K$.
Then the following are equivalent:
\begin{enumerate}
\item Every element of $\O$ is a $p$-th power $\vpmod{p\O}$. 
\item $(p, \disc K) = 1$.
\end{enumerate}
\end{lem}

We now state some results regarding orders in an algebraic number field.
\vskip3mm

\begin{defn}\label{defn:order}
An order $R$ in an algebraic number field $K$ is a subring of $K$
containing $1$ and which is a $\Z$-submodule of $K$ of maximum rank, i.e. of rank $n$ = $\rm{deg}(K/\Q)$.
\end{defn}

One notes that $\O$ is an order of $K$ and $\O$ contains every order, 
hence $\O$ is called the maximal order of $K$. The discriminant of an order
$R$ is defined to be the discriminant of any $\Z$-basis of $R$.

\begin{lem}\label{lem:basis}(\cite{KK2}, Lemma $4$)
If $K$ is a number field of degree $n$, $\O$ the ring of integers of $K$ and $R$
an order in $K$, then there is a $\Z$-basis $\alpha_1, \alpha_2, \ldots, \alpha_n$
of $R$ such that $\alpha_i = f_i\theta_i$, $f_i \in \Z$, $f_i > 0$ and moreover
$f_1 | f_2 | \cdots |f_n$. In fact with bases of $\O$ and $R$ as defined above, 
index of R in $\O = f_1f_2\cdots f_n$ and $\disc R = (f_1f_2\cdots f_n)^2\disc K.$
\end{lem}

\begin{lem}\label{lem:porder}(\cite{KK2}, Lemma $6$)
Let $R$ be a commutative ring with unity. Let $p$ be a prime and $R/pR$ a finite ring.
Then the following are equivalent:
\begin{enumerate}
\item Every element of $R$ is a $p$-th power $\vpmod{pR}$.
\item $x \in R, x^p \in pR$ implies $x \in pR$.
\end{enumerate}
\end{lem}

We now state a lemma from \cite{KK2} in the form that we would be using.
\begin{lem}\label{lem:ppower}
Let $R$ be an order in an algebraic number field $K$.
Let $p$ be a prime. 
Then the conditions $(1)$ and $(2)$ in Lemma \ref{lem:porder}
above are equivalent to the condition $(p, \disc R) = 1$.
\end{lem}
\begin{proof}

The proof that condition $(2)$ is equivalent to $(p, \disc R) = 1$
is a part of the proof of Theorem $1$ in \cite{KK2}. (See
\cite{KK2}, \S 3 for the details).
\end{proof}

\section{$2 \t 2$ matrices over orders as sums of cubes.}

\begin{lem}\label{lem:Three}
Let $R$ be a ring.
Let $A$ belong to $M_2(R)$. If $A = B^3$ then
there exist $t, \Delta$ in $R$ such that
$\trace A = t^3 - 3t\Delta$.
Conversely, given an element of $R$ of the form
$t^3 - 3t\Delta$, for some
$t, \Delta$ in $R$, there exists $B$ in $M_2(R)$
such that $\trace B^3 = t^3 - 3t\Delta$.
\end{lem}
\begin{proof}
Let $A =B^3$, $\trace B =t$ and $\det B = \Delta$. 
Note that $B$ satisfies its characteristic polynomial:
$B^2 - tB + \Delta I =0$, where $I$ stands for the $2 \t 2$ identity matrix.
This gives that $B^3 = tB^2 - \Delta B = t(tB-\Delta I)- \Delta B =
(t^2 - \Delta)B -t\Delta I$. Hence,
$A = B^3= (t^2 - \Delta)B -t\Delta I$ and 
$\trace A = t(t^2 - \Delta) - 2t \Delta = t^3 -3t\Delta$.

Conversely, given an element of $R$ of the form
$t^3 - 3t\Delta$, for some $t, \Delta$ in $R$, then
$B$ =
$
\begin{pmatrix}
  t   & -\Delta \\
  1   & 0 \\
\end{pmatrix}
$
is a matrix such that $\trace B^3 = t^3 - 3t\Delta$.
\end{proof}

\begin{lem}\label{lem:Ringcube}
Let $R$ be a ring. If every matrix in $M_2(R)$ 
is a sum of cubes of matrices in $M_2(R)$, then
every element of $R$ is a cube $\vpmod{3R}$.
\end{lem}
\begin{proof}
The proof follows from Lemma \ref{lem:Ringqth}. 
\end{proof}

{\it Here on in this section $R$ denotes an order in an algebraic
number field $K$.}

In the following lemma we exploit the structure on the elements of $R/3!R$ that can be expressed as sums of traces of cubes of $2 \t 2$ matrices over $R$ $\vpmod{3!R}$. Such a consideration seems natural from Corollary \ref{cor:Fequiv}. Here $\tilde{a}$ stands for the image of an element $a$ in $R$ under the map $R \longrightarrow R/3!R$.
\begin{lem}\label{lem:S}
Let
$$
S = \bigg \{ \tilde{\alpha} \in R/3!R~|~ {\alpha} = \s1,
~~{\rm where} ~t_i,~~ \Delta_i \in R \bigg \}$$

Then $S$ is a finite monoid and moreover an abelian group
generated additively by 
${\tilde t_m}^3,~ 3\tilde{t'_j}$ where $t_m,~ t'_j$ vary over $R$.
\end{lem}
\begin{proof}
As $R$ is an order in an algebraic number field, $S$ is finite. 
Clearly $S$ is an abelian monoid. 
Now if $\tilde{\alpha}$ belongs to $S$, then 
$3!\tilde{\alpha} = \overline{0}$ and so 
$\tilde{\alpha}+ (3! -1)\tilde{\alpha} = \overline{0}$.  
This shows that $S$ is a group.

To complete the proof it is enough to prove that the elements 
${\tilde t_m}^3,~ 3\tilde{t'_j}$ where $t_m,~ t'_j$ vary over $R$
are generated by elements of the form $(\tilde{t_i}^3- 3\tilde{t_i}\tilde{\Delta_i})$. 

Substituting $\Delta_i~=~0$ in  
$({t_i}^3~-~ 3{t_i}{\Delta_i})$, 
we get that ${\tilde t_i}^3$ is in $S$. Again putting 
${\Delta_i}~=~1$ in the same expression we get that 
${\tilde t_i}^3 - 3 {\tilde t_i}$ is in $S$. As $S$ is a 
group $3{\tilde t_i}$ is in $S$. 
\end{proof}

In the lemma below we see how $S$ characterizes those $2 \t 2$ matrices
over $R$ that are sums of cubes of $2 \t 2$ matrices over $R$. 
\begin{lem}\label{lem:Smatrix}
$M$ in $M_2(R)$ is a sum of cubes of matrices over $R$ if and only if
$\tilde{\alpha}$ belongs to $S$, where $\alpha = \trace{M}$.
\end{lem}
\begin{proof}
Let $M = \sum_{i=1}^{n}{M_i}^3$. Then, $\trace M = \sum_{i=1}^{n}{\trace(M_i)}^3 = \sum_{i=1}^{n}({t_i}^3 - 3{t_i}{\Delta_i})$.  
Thus if $\trace M = \alpha$ then reducing the above equation modulo $3!R$, 
we see that $\tilde{\alpha}$ belongs to $S$. 

Conversely if $\tilde{\alpha}$ belongs to $S$ where $\alpha = \trace{M}$, 
then $\tilde{\alpha} = 
\sum_{i=1}^{j_{\alpha}}(\tilde{{t_i}^3} - 3{\tilde{t_i}}{\tilde{\Delta_i}}) = \sum_{i=1}^{j_{\alpha}}{\tilde{\alpha_i}}$, where 
$\alpha_i = \trace {M_i}^3$, for some $M_i$'s in $M_2(R)$ 
(by Lemma \ref{lem:Three}) which translates
as $\trace (M - \sum_{i=1}^{j_{\alpha}} {M_i}^3)$ belongs to $3!R$.
Now Theorem \ref{thm:Trace} shows that $M$ is a sum of cubes of matrices over $R$.
\end{proof}

We now give a necessary and sufficient condition for $2 \t 2$
matrices over orders in an algebraic number field to be sums of
cubes of matrices.

\begin{thm}\label{thm:Rcube}
Let $K$ be an algebraic number field and $R$ an order in $K$.
Then, every matrix in $M_2(R)$
is a sum of cubes of matrices in $M_2(R)$ if and only if
$(3, \disc R) = 1$.
\end{thm}
\begin{proof}
$(\Longrightarrow)$
If every matrix in $M_2(R)$ is a sum of cubes of matrices in 
$M_2(R)$, then Lemma \ref{lem:Ringcube} shows that every
element of $R$ is a cube $\vpmod{3R}$. Now Lemma \ref{lem:ppower}
implies that $(3, \disc R) = 1$.

$(\Longleftarrow)$ Let $(3, \disc R) = 1$. Then Lemma \ref{lem:ppower} shows that every element of $R$ is a cube $\vpmod{3R}$. Thus,
in particular if $\trace M = \alpha$, then
$\alpha = \beta^3 + 3 \delta$, for some
$\beta$ and $\delta$ in $R$. On reducing modulo $3!R$ we get that
$\tilde\alpha$ belongs to $S$. This completes the proof using
Lemma \ref{lem:Smatrix}.
\end{proof}

\begin{rem}\label{rem:Zcube}
S. A. Katre and D. N. Sheth (see \cite{KS}) have obtained a criterion for
a $2 \t 2$ matrix over $R$, $R$ an integrally closed domain, to be 
a cube and have proved that every matrix in $M_2(\Z)$ is a sum of at 
most four cubes. For $M_2(\Z)$ whether $3$ cubes suffice is open.  
If $R$ is an order in an algebraic number field, 
and $(3, \disc R) = 1$ then the question of minimum number of 
cubes required to express a $2 \t 2$ matrix as a sum of cubes 
of matrices over $R$ is open. 
\end{rem}

\section{A counter-intuitive example}

Observe that (with the notations as defined in the introduction) if $\O_3 = \O$, then Corollary \ref{cor:T3} shows that every matrix in $M_2{\O}$ 
is a sum of cubes of matrices in $M_2{\O}$. That the converse does not
hold is seen from the example below. 
 
Note that Theorem \ref{thm:Stemm} stated in the next section gives
some conditions under which $\O_3 \neq \O$. In the particular
case of quadratic number fields, the conditions are as follows: $3$ is unramified and $2$ has a prime factor of degree $2$ in $\O$. It can be checked that these conditions hold for the quadratic
number field $K$ = $\mathbb Q(\sqrt5)$.

In this example let $K$ = $\mathbb Q(\sqrt5)$.  
Then $\O~=~\Z[\frac{1~+~\sqrt5}{2}]$ and $\disc K = 5$.

Let $\alpha = \frac{1 + \sqrt5}{2}$. Then, $\alpha$ satisfies the
recurrence relation $\alpha^2 - \alpha -1 = 0$.  
Thus, $\alpha^2 = \alpha + 1$ and $\alpha^3 = 2 \alpha + 1$.
Now if $\beta = a + b\alpha$ is an element of $\O$, then 
$\beta^3 = (a^3 + 3a^2b + b^3) + \alpha~(3ab(a+b) + 2b^3)$
which shows that the coefficient of $\alpha$ in the 
cube of any arbitrary element of $\O$ is always even.
This observation now shows that $\alpha$ cannot be 
expressed as a sum of cubes of elements of $\O$: for if 
$\alpha = {\sum_{j = 1}^{s}\beta_j^3}$, where $\beta_j$'s
belong to $\O$, then comparing the coefficient of $\alpha$ 
on both the sides of this equation, we get that $1$ is an 
even number, a contradiction. 
Thus $\alpha$ cannot be written as a sum of cubes of elements
of $\O$. However, Theorem \ref{thm:Rcube} shows that every $2 \t 2$ matrix 
over $\O$ is a sum of cubes of matrices over $\O$, as $(3, 5) = 1$.

\section{$2 \t 2$ matrices over orders as sums of fourth powers.}

\begin{lem}\label{lem:Four}
Let $R$ be a ring.
Let $A$ belong to $M_2(R)$. If $A = B^4$ then
there exist $t, \Delta$ in $R$ such that
$\trace A = t^4 - 4t^2\Delta + 2{\Delta}^2$.
Conversely, given an element of $R$ of the form
$t^4 - 4t^2\Delta + 2{\Delta}^2$, for some 
$t, \Delta$ in $R$, there exists $B$ in $M_2(R)$
such that $\trace B^4 = t^4 - 4t^2\Delta + 2{\Delta}^2$.  
\end{lem}
\begin{proof}
Let $A =B^4$ and $\trace B =t$, $\det B = \Delta$, then 
Note that $B$ satisfies its characteristic polynomial: 
$B^2 - tB + \Delta I =0$, where $I$ stands for the $2 \t 2$ identity matrix. This gives that $B^3 = tB^2 - \Delta B = t(tB-\Delta I)- \Delta B = 
(t^2 - \Delta)B -t\Delta I$. Hence, 
$A = B^4 = (t^2 - \Delta)B^2- t\Delta B = 
(t^2 - \Delta)(tB - \Delta I) - t\Delta B
= [t(t^2 - \Delta) - t \Delta]B - \Delta (t^2 - \Delta)I$.
Hence 
$\trace A = [t(t^2 - \Delta) - t \Delta]t - 2\Delta (t^2 - \Delta)
= t^4 - 4t^2 \Delta + 2{\Delta}^2$. 

Conversely, if there exists an element of $R$ of the type 
$t^4 - 4t^2\Delta + 2{\Delta}^2$, for some $t, \Delta$ in $R$, 
then the matrix 
$B$ =  
$
\begin{pmatrix}
  t   & -\Delta \\
  1   & 0 \\
\end{pmatrix} 
$
is a matrix such that $\trace B^4 = t^4 - 4t^2\Delta + 2{\Delta}^2$.
\end{proof}

\begin{lem}\label{lem:Reven}
Let $R$ be a ring and $n \geq 2$ be an integer.
If a matrix in $M_n(R)$ 
is a sum of $2j$-th powers of matrices in $M_n(R)$, where 
$j \geq 1$ (and thus a sum of squares), then $\trace M$ is a square $\vpmod{2R}$. Further if every matrix in $M_n(R)$ is a sum of $2j$-th powers of matrices in $M_n(R)$, ($j \geq 1$) then every element of 
$R$ is a square $\vpmod{2R}$.
\end{lem}

\begin{proof}
The proof follows from Lemma \ref{lem:Ringqth}. 
\end{proof}

{\it Here on in this section, $R$ denotes an order in an algebraic 
number field $K$.}  

In the following lemma we exploit the structure on the elements of 
$R/4!R$ that can be expressed as sums of traces of fourth powers of 
$2 \t 2$ matrices over $R$ $\vpmod{4!R}$. 
Such a consideration seems natural from Corollary \ref{cor:Fequiv}.
Here $\tilde{a}$ stands for the image of an element 
$a$ in $R$ under the map $R \longrightarrow R/4!R$. 

\begin{lem}\label{lem:G}
Let
$$G = \bigg \{ \tilde{\alpha} \in R/4!R~|~ {\alpha} = \w4 
,~~{\rm where}~ t_i, \Delta_i \in R \bigg \}$$
Then $G$ is a finite monoid and moreover an abelian group 
generated additively by 
${\tilde t_i}^4,~ 4{\tilde u_j}^2 {\tilde v_j}~, 2{\tilde \Delta_l}^2$
where $t_i, u_j, v_j, \Delta_l$ vary over $R$.    
\end{lem}
   
\begin{proof}
As $R$ is an order in an algebraic number field, $G$ is finite. 
Clearly $G$ is an abelian monoid. 
Now if $\tilde{\alpha}$ belongs to $G$, then $4!\tilde{\alpha} = \overline{0}$.
i.e. $\tilde{\alpha}+ (4! -1)\tilde{\alpha} = \overline{0}$.  
This shows that $G$ is a group. 

To complete the proof it is enough to prove that the elements 
${\tilde t_i}^4, 4{\tilde u_j}^2{\tilde v_j}, 2{\tilde{\Delta_l^2}}$, 
where $t_i, u_j, v_j, \Delta_l$ belong to $R$ are generated by elements 
of the form 
${\tilde t_i}^4 - 4{\tilde t_i}^2{\tilde \Delta_i} + 2{\tilde{\Delta_i^2}}$. 

Putting ${\Delta_i}=0$, in ${t_i}^4 - 4{t_i}^2{\Delta_i} + 2{\Delta_i^2}$
we have $\tilde{t_i}^4$ is in $G$.
Putting ${t_i}=0$, again in ${t_i}^4 - 4{t_i}^2{\Delta_i} + 2{\Delta_i^2}$, 
we get that $2\tilde{\Delta_i}^2$ is in $G$. As $G$ is a group, 
$4~\tilde{t_i}^2~{\tilde \Delta_i}$ belongs to $G$. 
\end{proof}

In the lemma below we see how $G$ characterizes those $2 \t 2$ matrices
over $R$ that are sums of fourth powers of $2 \t 2$ matrices over $R$. 

\begin{lem}\label{lem:Gmatrix}
$M$ in $M_2(R)$ is a sum of fourth powers of matrices over $R$ if and only if $\tilde{\alpha}$ belongs to $G$, where $\alpha = \trace{M}$.
\end{lem}
\begin{proof}
Let $M = \sum_{i=1}^{n}{M_i}^4$. Then, $\trace M = \sum_{i=1}^{n}{\trace({M_i}^4)} = \sum_{i=1}^{n}({t_i}^4~-~4~{t_i}^2~{\Delta_i}~+~2{\Delta_i}^2)$.  
If $\trace M = \alpha$, reducing the above equation modulo $4!R$,
shows that $\tilde{\alpha}$ belongs to $G$. 

Conversely if $\tilde{\alpha}$ belongs to $G$, where $\alpha = \trace M$,
then $\alpha = \sum_{i=1}^{j_{\alpha}}({t_i}^4 - 4{t_i}^2{\Delta_i} + 2{\Delta_i^2})$,
for $t_i, \Delta_i$ in $R$ by Lemma \ref{lem:Four}.  
So, $\tilde{\alpha} = \sum_{i=1}^{j_{\alpha}}({\tilde{t_i}}^4 - 4{\tilde{t_i}}^2{\tilde{\Delta}_i} + 2{\tilde{\Delta_i}^2})$
i.e. $\tilde{\alpha}=\sum_{i=1}^{j_{\alpha}}{\tilde{\alpha_i}}$, where 
$\alpha_i = \trace {M_i}^4$, for some $M_i$'s in $M_2(R)$ 
(by Lemma \ref{lem:Four}) which translates 
as $\trace (M - \sum_{i=1}^{j_{\alpha}} {M_i}^4)$ belongs to $4!R$. 
Now Theorem \ref{thm:Trace} shows that $M$ is a sum of fourth powers 
of matrices over $R$.  
\end{proof}

We now prove a necessary and sufficient condition for $2 \t 2$ matrices over orders to be sums of fourth powers of matrices. 

\begin{thm}\label{thm:Rfourth}
Let $K$ be an algebraic number field and $R$ an order of $K$.
Then, every matrix in $M_2(R)$ 
is a sum of fourth powers of matrices over $R$ if and only if 
$(2, \disc R) = 1$ i.e. discriminant of the order is odd. 
\end{thm}
\begin{proof}
$(\Longrightarrow)$
Let $R$ be an order of $K$. 
If every matrix in $M_2(R)$ is a sum of fourth powers of matrices in 
$M_2(R)$, then Lemma \ref{lem:Reven} shows that every
element of $R$ is a square $\vpmod{2R}$. Now Lemma \ref{lem:ppower}
implies that $(2, \disc R) = 1$.

$(\Longleftarrow)$ Let $(2, \disc R) = 1$. Then Lemma \ref{lem:ppower} shows that every element of $R$ is a square $\vpmod{2R}$. Thus, in particular if $\trace M = \alpha$, then $\alpha = \gamma^2 + 2 \delta$, 
for some $\gamma$ and $\delta$ in $R$. 
Repeated application of Lemma \ref{lem:ppower}
shows that $\alpha = ({\beta}^2 + 2\nu)^2 + 2({\theta}^2 + 2\lambda)$, 
for $\beta, \nu, \theta, \lambda \in R$ such that $\gamma = {\beta}^2 + 2\nu$ and $\delta = {\theta}^2 + 2\lambda$.
Thus $\alpha = {\beta}^4 + 4{\beta}^2\nu + 4{\nu}^2 
+ 2{\theta}^2 + 4\lambda$. Reducing modulo $4!R$ and applying 
Lemma \ref{lem:G} shows that $\tilde\alpha \in G$ and 
hence by Lemma \ref{lem:Gmatrix} it follows that every matrix in
$M_2(R)$ is a sum of fourth powers of matrices in $M_2(R)$. 
\end{proof}

\section{Matrices over the ring of integers of a quadratic number field
as sums of prime powers.}
We recall a theorem below:

\begin{thm}\label{thm:Stemm}(P. T. Bateman, R. M. Stemmler, \cite{BS})
Let $\O$ be the ring of integers of any algebraic number field
$K$ (not necessarily quadratic) and $q$ be a prime. Let
$\mathcal O_{q}$ be the additive subgroup generated by the $q$-th
powers of the elements of $\O$. Then, $\mathcal O_{q} \neq \O$ if
and only if at least one of the following holds:
\begin{enumerate}
\item $q$ is ramified in $\O$ (i.e. $q$ divides $\disc K$), 
\item $q$ is expressible in the form $\frac{p^r~ -~1}{p^d~-~1}$, 
where $p$ is a prime, $r$ and $d$ are positive integers, $d$ is a divisor of $r$  
and $p$ has in $\O$ a prime factor of degree $r$. 
\end{enumerate}
\end{thm}

We now prove the following theorem for the ring $\O$ of integers  
of a {\it quadratic} number field $K$:
\begin{thm}\label{thm:Oquad}
Let $n,q \geq 2$ and $q$ be a prime such that $n < q$.
Let $K$ be a quadratic number field and $\O$ its ring of 
integers. Then $(q, \disc K) = 1$ if and only if every 
$n \t n$ matrix over $\O$ 
is a sum of $q$-th powers of matrices over $\O$.
\end{thm}
\begin{proof}
If $q$ is unramified in $\O$
i.e. $(q, \disc K)=1$, then the only way in which 
${\O}_q \neq \O$ is when primes $q$ of the form listed in
Theorem \ref{thm:Stemm} occur.
However specialization to quadratic number fields shows that
there are restrictions on the integers $r$ and $d$ occurring in the 
expression for $q$.
In fact the only non-trivial case occurs when 
$r = 2$ and $d = 1$, in which case $q = \frac{p^2~-~1}{p~-~1} = p +1$.
This case occurs only when $q = 3$ and $p = 2$ and the additional 
assumption $n < q$ implies $n = 2$. For this case i.e. $n=2$ and $q=3$, the required result now follows from Theorem \ref{thm:Rcube}. 

Now the table in \cite{BS} shows that the primes $q$ of the above form 
are very rare. In fact, all such primes $q \geq 5$ are never of the form 
$\frac{p^2~-~1}{p~-~1}$, for some  prime $p$, having a prime factor
of degree $2$ in $\O$. All the exponents of $p$ in such prime powers 
are $\geq~3$. Thus, we have that for all primes $q \geq 5$, 
the only condition for $\O = \O_q$  
is that $(q, \disc K)= 1$. This observation along with Corollary \ref{cor:T3} proves that every $n \t n$ matrix over $\O$ is a sum of $q$-th powers of matrices over $\O$.

Conversely, if every matrix over $\O$ 
is a sum of $q$-th powers of matrices over $\O$, then 
in particular taking $R$ to be the ring $\O$ of integers of a 
quadratic number field $K$ in Lemma \ref{lem:Ringqth}, we have that
every element of $\O$ is a $q$-th power $\vpmod{q\O}$. 
Now Lemma \ref{lem:Disc} applied with $R = \O$ 
shows that $(q, \disc K) = 1$. 
\end{proof}

\begin{rem}\label{rem:Bound}
Let $n,q \geq 2$ be integers and $q$ be a prime. Let $K$ be an algebraic number field and $\O$ its ring of integers. 
Let $M \in M_n(\O)$ be a sum of $q$-th powers of matrices in $M_n(\O)$.
A rough bound on the number of $q$-th powers required to write $M$
as a sum of $q$-th powers may be given and it is five, if $n \ge q$
(see Theorem \ref{thm:Richman}) and
$g_k + 2^{q-1} + \frac{q-1}{3} + 1$ if $n < q$.
(See Theorem \ref{thm:Richman}, Theorem \ref{thm:Trace} and the results from \cite{S} viz. Theorem $9$ and Theorem $10$ for details). 
\end{rem}

Combining the results of Theorem \ref{thm:Oquad} and 
Remark \ref{rem:Bound} above, we have the following result for
{\it quadratic} number fields.
 
\begin{thm}\label{thm:together}
Let $\O$ be the ring of integers of a quadratic number field $K$.
Let $n,q \geq 2$ be integers and let $q$ be a prime. 
For $g_k$ as defined in Theorem \ref{thm:Trace}, we have: 
\begin{enumerate}
\item If $n \geq q$, then every $n \t n$ matrix is a sum of 
at most five $q$-th powers (by Theorem \ref{thm:Richman})
if and only if $(q, \disc K) = 1$.
\item If $n < q$ (so $q \geq 3$) then every $n \t n$ matrix over 
$\O$ is a sum of at most $g_k + 2^{q-1} + \frac{q-1}{3} + 1$ 
$q$-th powers of matrices over $\O$ if and only if 
$(q, \disc K) = 1$.
\end{enumerate}
\end{thm}
 
We give here some examples to illustrate some theorems that
have been proved in this paper.

\begin{enumerate}
\item The prime $11$ cannot be written in the form 
$\frac{p^r~-1}{p^d~-1}$ for any prime number $p$,
where $r$ and $d$ are positive integers and $d$ divides $r$. 
Therefore if $K$ is an algebraic number field (with 
ring of integers $\O$) such that 
$(11, \disc K) = 1$, Theorem \ref{thm:Stemm} 
shows that $\O_{11} = \O$. 
Hence by Corollary \ref{cor:T3} every $n \t n $ matrix 
(with $n \geq 2$, an integer) over $\O$ can be written 
as a sum of $11$-th powers matrices in $M_n{\O}$.

\item Let $K$ = $\Q(\sqrt[3]{2})$. In this case it can be checked that 
$\O = \Z[\sqrt[3]{2}]$. (See \cite{M}, Chapter 2, Exercise 41 or 
use the free mathematical software KASH (See \cite{kash}).
Now, $3 = \frac{2^2~-1}{2~-1}$. However $2$ does not have a prime factor
of degree $2$ in $\O$; since $2\O = (\alpha)^3$, for $\alpha = \sqrt[3]{2}$ 
in $\O$ and thus $(\alpha)$ is a prime ideal of $\O$ 
having residue degree $1$.
Thus in this case for every $2 \t 2$ matrix to be a sum of 
cubes of matrices over $\O$, it is required to see if  
$(3, \disc K) = 1$. Also using KASH, for $K = \Q(\sqrt[3]{2})$, 
we have $\disc K = -108$, and as $3$ divides $\disc K$, Theorem
\ref{thm:Rcube} shows that there are $2 \t 2$ matrices over $\O$ that 
cannot be written as sums of cubes of matrices over $\O$.

\item Let $K$ = $\Q(\sqrt[3]{19}, \omega)$ where $\omega = e^{2\pi i/3}$.
Let $\O$ denote its ring of integers. 
One has $13 = \frac{3^3~-1}{3~-1}$ and we show that 
$3$ cannot have a prime factor of degree $3$ in $\O$, so that one 
can appeal to Theorem \ref{thm:Stemm} if one also has that 
$(13, \disc K) = 1$. 

Note that $3\O$ is a product of three prime ideals each having 
ramification index $2$ and residue degree $1$ (See \cite{M}, 
Chapter 4, Page 103 or Factor $3*\O$ using KASH) and thus  
$3$ does not have a prime factor of degree $3$ in $\O$.
Also, $\disc K$ turns out to be
$(-3)^3 \cdot (-361)^2$, which is computed using the software
KASH and the transitivity property of discriminants
(see \cite{Ser}, Chapter 3, \S 4).
Thus $\disc K$ is coprime to $13$ and so
Theorem \ref{thm:Stemm} shows that $\O_{13} = \O$.
Now Corollary \ref{cor:T3} shows that every $n \t n $ matrix with
$n \geq 2$ can be written as a sum of $13$-th powers matrices
in $M_n{\O}$.

\item Let $K$ = $\Q(\sqrt[3]{7}, \sqrt[3]{11})$ and 
let $\O$ denote the ring of integers of $K$. 
Here $K_{0} = \Q(\sqrt[3]{11})$ 
and $K_{1} = \Q(\sqrt[3]{7})$ are of class numbers $2$ and 
$3$ respectively. Using KASH, one check that
$ \rm{Norm}_{\Q}(\rm{\disc {K_1({\sqrt[3]{11}})}/ K_1})  =  \left\langle 5314683 \right\rangle$ and
$\rm {\disc}(K_1/ \Q) = -1323$. Since $\left\langle\disc(K/\Q)\right\rangle=\left\langle\disc(K_0/\Q)\right\rangle^{[K:K_1]}\cdot{\rm Norm}_{\Q}(\disc(K/K_1))$ we get
$\left\langle \rm{\disc}(K/ \Q)\right\rangle =
{\left\langle -1323 \right\rangle}^3 \cdot \left\langle 5314683 \right\rangle$.
Since $K /\Q$ has one real and eight complex embeddings in $\C$,
sign of $\rm{\disc}(K/ \Q)$ is given by $(-1)^{8/2} = 1$, so 
$\rm{\disc}(K/ \Q) = {1323^3}\cdot 5314683$. Thus 
$\rm{\disc}(K/ \Q)= (3)^{10}\cdot(7)^{6}\cdot(11)^6$.
As $(2, \disc K) =1$, Theorem \ref{thm:Rfourth} proves that
every $2 \t 2$ matrix over $\O$ is a sum of fourth powers
(and a sum of squares, too) of matrices over $\O$.
Also as $3$ divides $\disc K$, Theorem \ref{thm:Rcube} shows 
that are $2 \t 2$ matrices over $\O$ that are not sums of cubes 
of matrices over $\O$.
\end{enumerate}

\medskip
\it{Acknowledgement}: 
We thank R. Balasubramanian for pointing to us some work of
mathematicians on Waring's problem for algebraic number fields
and Dinesh Thakur for informing us about the free
mathematical software KASH, useful for computations in algebraic
number theory. We also thank Bhaskaracharya Pratishthana, Pune for some facilities. The first author would like to thank the Council of Scientific and Industrial Research (C.S.I.R), New Delhi for awarding a research fellowship during her doctoral 
studies.

\end{document}